\documentstyle[amscd,amssymb,verbatim,epsf]{amsart}

\begin{document}
\theoremstyle{plain}
\newtheorem{Thm}{Theorem}
\newtheorem{Cor}{Corollary}
\newtheorem{Ex}{Example}
\newtheorem{Con}{Conjecture}
\newtheorem{Main}{Main Theorem}
\newtheorem{Lem}{Lemma}
\newtheorem{Prop}{Proposition}

\theoremstyle{definition}
\newtheorem{Def}{Definition}
\newtheorem{Note}{Note}

\theoremstyle{remark}
\newtheorem{notation}{Notation}
\renewcommand{\thenotation}{}

\errorcontextlines=0
\numberwithin{equation}{section}
\renewcommand{\rm}{\normalshape}%

\title[Geodesic Flow]%
   {Geodesic Flow on Global Holomorphic Sections of ${\mbox{TS}}^2$}
\author{Brendan Guilfoyle}
\address{Brendan Guilfoyle\\
          Department of Mathematics and Computing \\
          Institute of Technology, Tralee \\
          Clash \\
          Tralee  \\
          Co. Kerry \\
          Ireland.}
\email{brendan.guilfoyle@@ittralee.ie}
\author{Wilhelm Klingenberg}
\address{Wilhelm Klingenberg\\
 Department of Mathematical Sciences\\
 University of Durham\\
 Durham DH1 3LE\\
 United Kingdom.}
\email{wilhelm.klingenberg@@durham.ac.uk }

\date{17th February, 2006}

\begin{abstract}
We study the geodesic flow on the global holomorphic sections of the bundle
$\pi:{\mbox{TS}}^2\rightarrow \mbox{S}^2$ induced by the neutral K\"ahler metric on the
space of oriented lines of ${\Bbb{R}}^3$, which we identify with ${\mbox{TS}}^2$. This flow is shown to be
completely integrable when the sections are symplectic and the behaviour of the geodesics is 
described.
\end{abstract}

\maketitle

\section{Introduction}

It has long been known that certain global holomorphic sections of the bundle
$\pi:{\mbox{TS}}^2\rightarrow \mbox{S}^2$ correspond to the set of oriented lines through a point in 
${\Bbb{R}}^3$ \cite{hitch}.
However, this consists of only a three real parameter family of what is a six parameter family of 
global holomorphic sections. The purpose of this paper is to consider the other three parameter family of
``twisting'' holomorphic spheres and their geometry. In particular, we will show that these spheres, which correspond
to overtwisted contact structures, exhibit interesting dynamical properties (see for example \cite{arn} 
\cite{moser} \cite{tab}). 

More precisely, we study the geodesic flow on the global holomorphic sections of the bundle
$\pi:{\mbox{TS}}^2\rightarrow \mbox{S}^2$ induced by the neutral K\"ahler structure on ${\mbox{TS}}^2$. 
Here, ${\mbox{TS}}^2$ is identified with the space ${\Bbb{L}}$ of oriented affine lines in ${\Bbb{R}}^3$ and
the K\"ahler metric is invariant under the induced action on ${\Bbb{L}}$ of the Euclidean group \cite{gak4}. 

For topological reasons such sections are quadratic and, modulo the Euclidean action, we reduce this to
a 1-parameter family. Within this family lies the oriented normal line congruence to the
round sphere in ${\Bbb{R}}^3$. This is the unique lagrangian surface within the class, the other line 
congruences being symplectic - having twist, and hence we refer to them as twisting holomorphic spheres.

The metric induced on the lagrangian surface is degenerate at each point and hence we study the twisting case.
We show that the lagrangian points on a twisting global holomorphic section project by the bundle map to a great 
circle on ${\mbox S}^2$. The induced metric is positive and negative definite on either side of this circle of lagrangian 
points, and degenerate on the circle. We study the geodesic flow on these two discs and prove that:

\vspace{0.2in}
\noindent{\bf Main Theorem}:

{\it The geodesic flow on twisting holomorphic spheres is completely integrable. Geodesics with zero 
angular momentum reach the lagrangian circle within a finite time, where the flow blows up, while geodesics with
non-zero angular momentum oscillate between a maximum and minimum distance from the circle.}

\vspace{0.2in} 

In the next section we give a summary of the geometry of the space of oriented affine lines in ${\Bbb{R}}^3$ - further
details can be found in \cite{gak4}. 
The following section describes the space of global holomorphic sections, while we prove the main result about the 
geodesic flow in Section 4.

\section{The Neutral K\"ahler Metric}

The space ${\Bbb{L}}$ of oriented lines in ${\Bbb{R}}^3$ can be identified with the tangent bundle to the 2-sphere
\cite{hitch}. This identification gives a useful way to compute the wealth of geometric structure that
exists on ${\Bbb{L}}$. In particular, using the complex coordinate $\xi$ on ${\mbox S}^2-\{\mbox{south pole}\}$ obtained 
by stereographic projection from the south pole, one obtains local complex coordinates ($\xi,\eta$) on ${\Bbb{L}}$ by 
the identification:
\[
(\xi,\eta)\leftrightarrow \eta\frac{\partial}{\partial \xi}+\bar{\eta}\frac{\partial}{\partial \bar{\xi}}\in
   {\mbox T}_\xi \;{\mbox S}^2.
\]

The Euclidean group acting on ${\Bbb{R}}^3$ sends oriented lines to oriented lines, and therefore induces an
action on ${\Bbb{L}}$. In our coordinates, a translation acts by:
\begin{equation}\label{e:trans}
\xi\rightarrow\xi'=\xi  \qquad 
\eta\rightarrow \eta'=\eta+\alpha_1-a_1\xi-\bar{\alpha}_1\xi^2,
\end{equation}
for $\alpha_1\in{\Bbb{C}},a_1\in{\Bbb{R}}$, while a rotation acts by
\begin{equation}\label{e:rot}
\xi\rightarrow\xi'=\frac{\alpha_2\xi+\alpha_3}{-\bar{\alpha}_3\xi+\bar{\alpha}_2}  \qquad 
\eta\frac{\partial}{\partial \xi}\rightarrow \eta'\frac{\partial}{\partial \xi'}=\frac{\eta}{(-\bar{\alpha}_3\xi+\bar{\alpha}_2)^2}\frac{\partial}{\partial \xi'},
\end{equation}
for $\alpha_2,\alpha_3\in{\Bbb{C}}$ satisfying $\alpha_2\bar{\alpha}_2+\alpha_3\bar{\alpha}_3=1$.

The natural geometric structure on ${\Bbb{L}}$ includes the projection $\pi:{\Bbb{L}}\rightarrow {\mbox S}^2$ which sends
an oriented line to its direction. We will be interested in global sections of this bundle, that is maps 
$s:{\mbox S}^2\rightarrow {\Bbb{L}}$ such that $\pi\circ s$ is the identity on ${\mbox S}^2$.

We also have a complex structure: ${\Bbb{J}}:{\mbox T}{\Bbb{L}}\rightarrow {\mbox T}{\Bbb{L}}$ such that 
${\Bbb{J}}\circ{\Bbb{J}}=-{\mbox{Id}}$ (plus an integrability condition).
 This is compatible with the above complex coordinates in the sense that
\[
{\Bbb{J}}\left(\frac{\partial}{\partial
    \xi}\right)=i\frac{\partial}{\partial \xi} \qquad \qquad
  {\Bbb{J}}\left(\frac{\partial}{\partial
    \eta}\right)=i\frac{\partial}{\partial \eta}.
\]

In addition, there is a natural symplectic structure on ${\Bbb{L}}$. That is, there is a non-degenerate closed 2-form on
${\Bbb{L}}$, which, in our coordinates, has the local expression:
\begin{equation}
\Omega=\frac{2}{(1+\xi\bar{\xi})^2}\left(
  d\eta \wedge d\bar{\xi}+d\bar{\eta} \wedge d\xi
   +\frac{2(\xi\bar{\eta}-\bar{\xi}\eta)}{1+\xi\bar{\xi}}d\xi\wedge d\bar{\xi}
\right).
\end{equation}

Finally, there exists a canonical K\"ahler metric ${\Bbb{G}}$ on ${\Bbb{L}}$ which is compatible with this complex
structure and is invariant under the action induced on ${\Bbb{L}}$ by the Euclidean isometry group acting on ${\Bbb{R}}^3$.
This has local expression:

\begin{equation}\label{e:metric}
{\Bbb{G}}=\frac{2i}{(1+\xi\bar{\xi})^2}\left(
  d\eta \otimes d\bar{\xi}-d\bar{\eta}\otimes d\xi
   +\frac{2(\xi\bar{\eta}-\bar{\xi}\eta)}{1+\xi\bar{\xi}}d\xi\otimes d\bar{\xi}
\right).
\end{equation}

\begin{Def}
A 2-parameter family of oriented lines in ${\Bbb{R}}^3$ forms a surface $\Sigma$ in ${\Bbb{L}}$ - which we refer to
as a {\it line congruence}. 
\end{Def}

We are interested in the geometric structures induced on $\Sigma$ by the K\"ahler 
structure (${\Bbb{J}}, \Omega, {\Bbb{G}}$). For the symplectic structure, we have:

\begin{Thm} \cite{arn}

A line congruence $\Sigma\subset{\Bbb{L}}$ is lagrangian (i.e. $\Omega|_\Sigma=0$) iff there exist surfaces
in ${\Bbb{R}}^3$ orthogonal to the lines.
\end{Thm}

Hence, we refer to a non-lagrangian (or symplectic) line congruence as {\it twisting}.
On the other hand, a line congruence $\Sigma$ is said to holomorphic if ${\Bbb{J}}$ preserves the tangent space 
of $\Sigma$. A holomorphic line congruence that is the graph of a local section can be described by a
holomorphic equation $\eta=\eta(\xi)$.

The signature of the metric induced on $\Sigma$ can be positive definite, negative definite, lorentz or degenerate.
In particular, we have the following:

\begin{Thm} \cite{gak4}

The metric induced on a holomorphic line congruence is either positive or negative definite or degenerate. It is 
degenerate precisely at the points on the line congruence where the symplectic form vanishes.
\end{Thm}

Generically, the lagrangian points on a holomorphic surface form curves, the only holomorphic line congruences
which are lagrangian everywhere being the oriented normals to planes and spheres in ${\Bbb{R}}^3$.

\section{Twisting Holomorphic Spheres}

Let us consider global holomorphic sections of the bundle $\pi:{\Bbb{L}}\rightarrow \mbox{S}^2$. Since the 
tangent bundle to the 2-sphere is of degree 2, such sections are quadratic:
\[
\eta=\beta_1+\beta_2\xi+\beta_3\xi^2,
\] 
for $\beta_1,\beta_2,\beta_3\in{\Bbb{C}}$. The Euclidean action on ${\Bbb{L}}$ allows us to put such sections into 
the standard form:

\begin{Prop}
After a rotation and a translation quadratic holomorphic sections can be put in the form:
\begin{equation}\label{e:holsph}
\eta=ci\xi,
\end{equation}
for $c\in[0,\infty)$.
\end{Prop} 
\begin{pf}

By a translation (\ref{e:trans}) with 
$\alpha_1={\scriptstyle{\frac{1}{2}}}(\bar{\beta}_3-\beta_1)$ and $a_1={\scriptstyle{\frac{1}{2}}}(\bar{\beta}_2+\beta_2)$
we can reduce this to
\[
\eta=\gamma+ci\xi+\bar{\gamma}\xi^2,
\] 
for $\gamma\in{\Bbb{C}}$ and $c\in{\Bbb{R}}$.
Now consider the rotation ({\it cf.} equation (\ref{e:rot})) with 
\[
\alpha_2=\frac{1}{(1+\xi_0\bar{\xi}_0)^{\scriptstyle{\frac{1}{2}}}} \qquad\qquad
\alpha_3=-\frac{\xi_0}{(1+\xi_0\bar{\xi}_0)^{\scriptstyle{\frac{1}{2}}}},
\]
for
\[
\xi_0=\frac{c-\left(c^2+4\gamma\bar{\gamma}\right)^{\scriptstyle{\frac{1}{2}}}}{2i\bar{\gamma}}.
\]
This rotation induces the change:
\[
\eta\rightarrow\eta'=\left(c^2+4\gamma\bar{\gamma}\right)^{\scriptstyle{\frac{1}{2}}}i\xi'.
\]
A relabelling of coordinates and constant yields the stated result.
\end{pf}

We now identify the lagrangian spheres in this class:

\begin{Prop}
A holomorphic section of the form (\ref{e:holsph}) is lagrangian iff $c=0$.
\end{Prop}
\begin{pf}
A line congruence $\Sigma\subset{\Bbb{L}}$ is lagrangian iff the symplectic form pulled back to $\Sigma$ vanishes.
For a section of the form (\ref{e:holsph}), we get
\[
\Omega|_\Sigma=\frac{2}{(1+\xi\bar{\xi})^2}\left[2ci-4ci\frac{\xi\bar{\xi}}{1+\xi\bar{\xi}}\right]d\xi\wedge d\bar{\xi}
=\frac{4ci(1-\xi\bar{\xi})}{(1+\xi\bar{\xi})^3}d\xi\wedge d\bar{\xi}.
\]
This vanishes on the whole of $\Sigma$ iff $c=0$.
\end{pf}

When $c=0$ the section represents the set of oriented lines through the origin in ${\Bbb{R}}^3$, or,
equivalently, the oriented normals to the round sphere centred at the origin. 

For $c\neq0$ the section is ``twisting'' in the sense there is no 
surface in ${\Bbb{R}}^3$ which is orthogonal to the lines. The line congruence can be viewed as follows: starting 
with the oriented line pointing along the positive $x^3$-axis, as one moves out from the axis, the line rotates 
until it is contained in the $x^1x^2$-plane (at a perpendicular distance $c$ from the origin). 
Then moving back towards the $x^3$-axis the line continues to rotate 
until, when it returns to the $x^3$-axis, it is pointing downwards (see the Figure below). 

\vspace{0.1in}
\setlength{\epsfxsize}{4.5in}
\begin{center}
   \mbox{\epsfbox{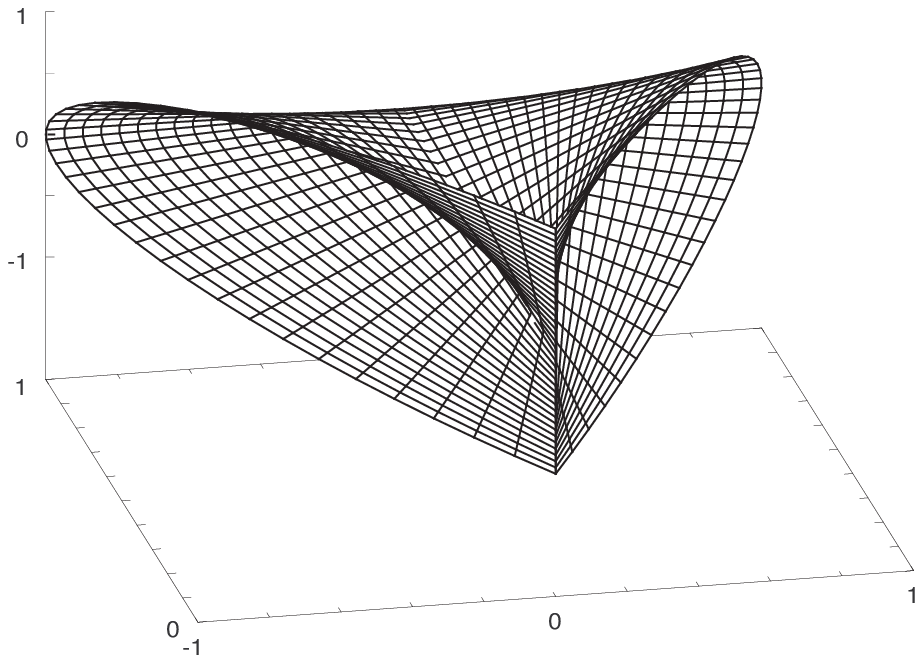}}
\end{center}
\vspace{0.1in}

This path on $\Sigma$ projects to a great circle from the north pole to the south pole on ${\mbox S}^2$. To see the full
2-parameter family of oriented lines we must rotate this about the $x^3$-axis. In fact, the planes orthogonal to the 
lines in this congruence form a distribution in ${\Bbb{R}}^3$ which is exactly that
of a pair of overtwisted contact structures \cite{eliah}.

From the above computation we also have that:

\begin{Prop}
The lagrangian points on the sphere $(\xi,\eta=ci\xi)$ for $c\neq0$ lie on the equator $|\xi|=1$.
\end{Prop}

\section{The Geodesic Flow}

We now look at the metric induced by the neutral K\"ahler metric on the twisting holomorphic spheres:

\begin{Prop}\label{p:fint}
The metric on the sphere $(\xi,\eta=ci\xi)$ is in local coordinates
\[
ds^2=-\frac{4c(1-\xi\bar{\xi})}{(1+\xi\bar{\xi})^3}d\xi\otimes d\bar{\xi}.
\]
Thus, for $c>0$ the metric is negative definite on the upper hemisphere, positive definite on the lower hemisphere and 
degenerate on the equator. 
\end{Prop}
\begin{pf}
This follows from pulling the metric (\ref{e:metric}) back to the line congruence.
\end{pf}

We turn now to the geodesic flow:

\begin{Prop}\label{p:main1}
Consider the holomorphic sphere $\Sigma\subset{\Bbb{L}}$ given by $(\xi,\eta=ci\xi)$ 
for $c>0$. The geodesic flow on $\Sigma$ is completely integrable with first integrals
\[
{\mbox {I}}_1=\frac{1-\xi\bar{\xi}}{(1+\xi\bar{\xi})^3}\;\dot{\xi}\dot{\bar{\xi}}
\qquad\qquad
{\mbox {I}}_2=\frac{1-\xi\bar{\xi}}{2i(1+\xi\bar{\xi})^3}\;\left(\bar{\xi}\dot{\xi}-\xi\dot{\bar{\xi}}\right).
\]
\end{Prop}
\begin{pf}
Consider the affinely parameterised geodesic $t\mapsto(\xi(t),ci\xi(t))$ on $\Sigma$ with tangent vector
\[
{\mbox T}=\dot{\xi}\frac{\partial}{\partial\xi}+\dot{\bar{\xi}}\frac{\partial}{\partial\bar{\xi}}.
\]
The geodesic equation 
${\mbox T}^j\nabla_j{\mbox T}^k=0$, projected onto the $\xi$ coordinate is
\[
\ddot{\xi}+\Gamma_{\xi\xi}^\xi\dot{\xi}^2 +2\Gamma_{\bar{\xi}\xi}^\xi\dot{\xi}\dot{\bar{\xi}}  
                       + \Gamma_{\bar{\xi}\bar{\xi}}^\xi\dot{\bar{\xi}}^2=0.
\] 
For the induced metric (as given in Proposition \ref{p:fint}) a straight-forward calculation yields the 
Christoffel symbols:
\[
\Gamma_{\xi\xi}^\xi=\partial\left[\ln\left(\frac{1-\xi\bar{\xi}}{(1+\xi\bar{\xi})^3}\right)\right]
\qquad\qquad
\Gamma_{\bar{\xi}\xi}^\xi=0 \qquad\qquad \Gamma_{\bar{\xi}\bar{\xi}}^\xi=0,
\]
where for short we have written $\partial$ for the derivative with respect to $\xi$. 
Thus the geodesic equation reduces to
\[
\ddot{\xi}=-\partial\left[\ln\left(\frac{1-\xi\bar{\xi}}{(1+\xi\bar{\xi})^3}\right)\right]\dot{\xi}^2.
\]

The fact that ${\mbox I}_1$ is constant along a geodesic comes from the fact that the geodesic flow preserves 
the length of the 
tangent vector ${\mbox T}^j$. On the other hand, differentiating ${\mbox I}_2$ with respect to $t$:
\begin{align}
2i\dot{{\mbox I}}_2=& \left[\partial\left(\frac{1-\xi\bar{\xi}}{(1+\xi\bar{\xi})^3}\right)\dot{\xi}
+\bar{\partial}\left(\frac{1-\xi\bar{\xi}}{(1+\xi\bar{\xi})^3}\right)\dot{\bar{\xi}}\right]
   \left(\bar{\xi}\dot{\xi}-\xi\dot{\bar{\xi}}\right)  
\nonumber\\
&\qquad\qquad+ \frac{1-\xi\bar{\xi}}{(1+\xi\bar{\xi})^3}
\left(\dot{\bar{\xi}}\dot{\xi}+\bar{\xi}\ddot{\xi}-\dot{\xi}\dot{\bar{\xi}}-\xi\ddot{\bar{\xi}}\right)   \nonumber\\
&=\partial\left(\frac{1-\xi\bar{\xi}}{(1+\xi\bar{\xi})^3}\right)
   \left(\bar{\xi}\dot{\xi}^2-\xi\dot{\xi}\dot{\bar{\xi}}\right) 
+\bar{\partial}\left(\frac{1-\xi\bar{\xi}}{(1+\xi\bar{\xi})^3}\right)
   \left(\xi\dot{\bar{\xi}}^2-\bar{\xi}\dot{\xi}\dot{\bar{\xi}}\right)  
\nonumber\\
&\qquad-\frac{1-\xi\bar{\xi}}{(1+\xi\bar{\xi})^3}
   \left\{\partial\left[\ln\left(\frac{1-\xi\bar{\xi}}{(1+\xi\bar{\xi})^3}\right)\right]\bar{\xi}\dot{\xi}^2
     -\bar{\partial}\left[\ln\left(\frac{1-\xi\bar{\xi}}{(1+\xi\bar{\xi})^3}\right)\right]\xi\dot{\bar{\xi}}^2\right\}             \nonumber\\
&=\frac{2\xi\bar{\xi}(2-\xi\bar{\xi})}{1-\xi^2\bar{\xi}^2}\dot{\xi}\dot{\bar{\xi}}- 
      \frac{2\xi\bar{\xi}(2-\xi\bar{\xi})}{1-\xi^2\bar{\xi}^2}\dot{\xi}\dot{\bar{\xi}} \nonumber\\
&=0, \nonumber
\end{align}
as claimed.
\end{pf}

It is clear that ${\mbox I}_2$ is angular momentum and the integrability of the
geodesic flow comes from conservation of this angular momentum.

The qualitative behaviour of the geodesic flow can now be determined:

\begin{Prop}\label{p:main2}
Geodesics with zero generalised
angular momentum ${\mbox I}_2$ reach the lagrangian circle within a 
finite time, where the flow blows up, while geodesics with
non-zero angular momentum oscillate between a maximum and minimum distance from the circle.
\end{Prop}
\begin{pf}
Let $\xi=Re^{i\theta}$, so that
\[
{\mbox {I}}_1=\frac{1-R^2}{(1+R^2)^3}\left(\dot{R}^2+R^2\dot{\theta}^2\right)
\qquad\qquad
{\mbox {I}}_2=\frac{1-R^2}{2i(1+R^2)^3}\;R^2\dot{\theta}.
\]
We will work in the upper hemisphere ($R<1$) so that ${\mbox I}_1\ge0$ - a similar analysis will hold in the lower hemisphere.

Let us first assume that the angular momentum is zero: ${\mbox I}_2=0$. Thus $\dot{\theta}=0$ throughout the motion and
so the geodesic is an arc of a great circle through the north pole. Integrating the first integral ${\mbox I}_1$ we get
\[
\sqrt{{\mbox I}_1}t+C_0=\int\frac{(1-R^2)^{\scriptstyle{\frac{1}{2}}}}{(1+R^2)^{\scriptstyle{\frac{3}{2}}}}dR.
\]
Now, the integral on the right hand side is a special case of the Appell hypergeometric function $f_1$:
\begin{align}
\sqrt{{\mbox I}_1}t+C_0=&R\;f_1(0.5;-0.5,1.5;1.5,R^2,-R^2) \nonumber\\
  &=\sum_{k=0}^\infty\sum_{l=0}^\infty(-1)^l\frac{(0.5)_{k+l}(-0.5)_k(1.5)_l}{(1.5)_{k+l}\;k!\;l!}R^{2(k+l)+1},\nonumber
\end{align}
where we have used the Pochhammer symbol $(a)_k=\Gamma(a+k)/\Gamma(a)$. Thus, starting from $R=0$, a geodesic reaches
the boundary $R=1$ in time $t\approx0.599070\;{\mbox I}_1^{\scriptstyle{-\frac{1}{2}}}$.

Now suppose that ${\mbox I}_2\neq0$. We then have that
\[
{\mbox I}_1-\frac{(1+R^2)^3}{(1-R^2)R^2}{\mbox I}_2^2=\frac{1-R^2}{(1+R^2)^3}\dot{R}^2.
\] 
Thus for $\dot{R}$ to remain real, we must have
\[
U_{\mbox eff}=\frac{(1+R^2)^3}{(1-R^2)R^2}\leq\frac{{\mbox I}_1}{{\mbox I}_2^2}.
\]
As the plot of $U_{\mbox eff}(R)$ below illustrates the geodesic oscillates between a maximum and minimum value for $R$.

\vspace{0.1in}
\setlength{\epsfxsize}{4.5in}
\begin{center}
   \mbox{\epsfbox{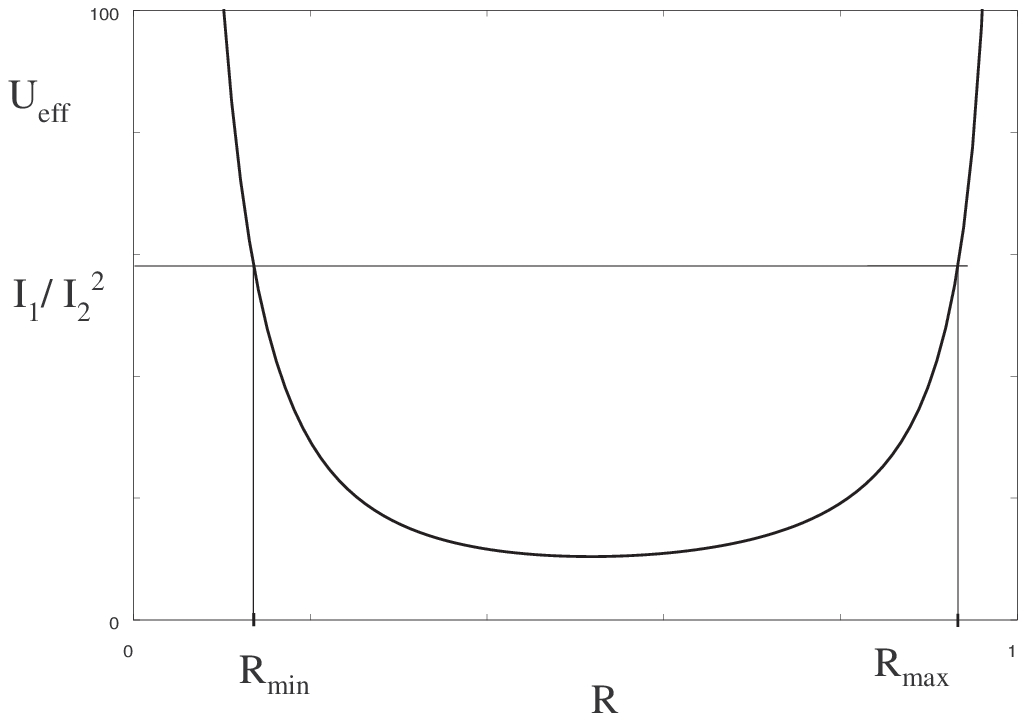}}
\end{center}
\vspace{0.1in}
\end{pf}

The main theorem follows from Propositions \ref{p:main1} and \ref{p:main2}.

\end{document}